\documentclass[11pt,reqno]{amsart}
\usepackage{amsfonts}
\usepackage{fancyhdr}
\usepackage{times}
\usepackage[T1]{fontenc}
\usepackage{mathrsfs}
\usepackage{latexsym}
\usepackage[dvips]{graphics}
\usepackage{epsfig}
\usepackage{hyperref, amsmath, amsthm, amsfonts, amscd, flafter,epsf}
\usepackage{amsmath,amsfonts,amsthm,amssymb,amscd}
\input amssym.def
\input amssym.tex
\usepackage{color}
\usepackage[parfill]{parskip}
 \usepackage{mathrsfs} % use command \mathscr
\usepackage{geometry}                % See geometry.pdf to learn the layout options. There are lots.
\usepackage{epstopdf}

\usepackage{verbatim}

\newtheorem{thm}{Theorem}%[section]
\newtheorem{cor}[thm]{Corollary}

\newtheorem{lem}[thm]{Lemma}
\newtheorem{thm*}{Theorem}
\newtheorem*{con*}{Conjecture}
\newtheorem*{lem*}{Lemma}
\newtheorem{prop}[thm]{Proposition}

\def\g{\gamma}

\begin{document}

\title{The distribution of zeros of $\zeta'(s)$ and gaps between zeros of $\zeta(s)$}
\author{Fan Ge}

\email{fan.ge@rochester.edu}

\address{Department of Mathematics, University of Rochester, Rochester, NY}

\maketitle

 \baselineskip=17pt
\hbox{}
\medskip

\begin{abstract}
Assume the Riemann Hypothesis, and let $\g^+>\g>0$ be ordinates of
two consecutive zeros of $\zeta(s)$. It is shown that if
$\g^+-\g<v/\log \g$ with $v<c$ for some absolute positive constant
$c$, then the box
$$
\{s=\sigma+it: 1/2<\sigma<1/2+v^2/4\log\g, \g\le t\le \g^+\}
$$
contains exactly one zero of $\zeta'(s)$. In particular, this allows
us to prove half of a conjecture of Radziwi{\l}{\l} \cite{Rad} in a
stronger form. Some related results on zeros of $\zeta(s)$ and
$\zeta'(s)$ are also obtained.
\end{abstract}

\section{Introduction}

%Let $N=N(T)$ denote the number of zeros of $\zeta(s)$ in the region
%$0< \sigma, T\le t\le 2T$, and let $N_1=N_1(T)$ be that of
%$\zeta'(s)$ in the same region.

Throughout this paper $s=\sigma+it$ is a complex variable. Write
$\rho=\beta+i\g$ and $\rho'=\beta'+i\g'$ to be a generic zero of
$\zeta(s)$ and $\zeta'(s)$, respectively. If $\zeta(1/2+i\g)=0$, let
$\g^+$ denote the smallest $t>\g$ with $\zeta(1/2+it)=0$. The phrase
"$\g$ is large" stands for "$\g$ is larger than an absolute
constant". Finally, we order the ordinates of zeros of $\zeta(s)$ as
$0<\g_1\le \g_2\le \cdots$, and similarly for zeros of $\zeta'(s)$.

The distribution of zeros of $\zeta'(s)$, and its relationship to
zeros of $\zeta(s)$, have been investigated by many authors (see
\cite{Ber}, \cite{CGh}, \cite{FKi}, \cite{Fen}, \cite{Guo},
\cite{Guo2}, \cite{GYi}, \cite{Ki}, \cite{LMo}, \cite{Rad},
\cite{Spe}, \cite{Sou}, \cite{Zha}). For example, a well-known
theorem of A. Speiser~\cite{Spe} states that the Riemann Hypothesis
(RH) is equivalent to $\zeta'(s)$ having no zeros in $0<\sigma<1/2$.

In \cite{Sou} K. Soundararajan raised the following conjecture.
\begin{con*}(Soundararajan.)
Assume RH. The following two statements are equivalent:
\begin{align*}&\textnormal{(A)}\quad
\liminf_{\g'\rightarrow\infty}(\beta'-1/2)\log\g'=0;\\
&\textnormal{(B)}\quad
\liminf_{\g\rightarrow\infty}(\g^+-\g)\log\g=0.
\end{align*}
\end{con*}

As stated by Soundararajan, both of these two assertions are almost
certainly true, and the point of this conjecture is that there might
be a simple way of relating (A) and (B) without reference to their
individual validity.

In \cite{Zha} Y. Zhang showed that on RH (B) implies (A), which
solved one direction of Soundararajan's Conjecture. By contrast, the
attempts to prove that (A) implies (B) have been unsuccessful (see
Theorem \ref{thm sound conj} below). As an alternative, D. W. Farmer
and H. Ki~\cite{FKi} considered the following. Let $w(x)$ be the
indicator function of the unit interval $[0, 1]$, and define
$$
m'(v)=\liminf_{T\rightarrow \infty}\frac{2\pi}{T\log T}\sum_{T\le
\g'\le 2T} w\bigg(\frac{(\beta'-1/2)\log T}{v}\bigg)
$$
and
$$
m(v)=\liminf_{T\rightarrow \infty}\frac{2\pi}{T\log T}\sum_{T\le
\g_n\le 2T} w\bigg(\frac{(\g_{n+1}-\g_n)\log T}{2\pi v}\bigg).
$$
They conjectured that if $m'(v)\gg v^\alpha$ for some $\alpha<2$,
then $m(v)>0$ for all $v>0$. This conjecture was viewed as a
refinement of Soundararajan's conjecture, and was recently proved
(on RH) by M. Radziwi{\l}{\l}~\cite{Rad}.

On the other hand, by investigating random matrix models for
$\zeta'(s)$, F. Mezzadri \cite{Mez} conjectured an asymptotic
formula for $m'(v)$. Namely,
\begin{align}\label{eq asymp m'}
m'(v)\sim \frac{8}{9\pi}v^{3/2}\ .
\end{align}
We also refer the readers to \cite{DFFHMP} for a detailed study in
this direction. The corresponding conjecture for $m(v)$ (see e.g.
\cite{FKi}) is
\begin{align}\label{eq asymp m}
m(v)\sim \frac{\pi}{6}v^3\ .
\end{align}
The validity of \eqref{eq asymp m'} and \eqref{eq asymp m} will have
significant consequences. In particular, from the work of J. B.
Conrey and H. Iwaniec \cite{CoI} and further work of Farmer and Ki
\cite{FKi}, either of \eqref{eq asymp m'} or \eqref{eq asymp m} (or
even weaker forms) will imply the non-existence of Landau-Siegel
zeros.

Our first result supports the relative size between $m(v)$ and
$m'(v)$ in the above conjectures. In particular, we prove one
direction of Radziwi{\l}{\l}'s conjecture, which asserts that if
$m(v)$ or $m'(v)$ is $\gg v^A$ then $m(v/2\pi)\asymp m'(v^2)$ (see
\cite{Rad}). (Note that we do not need the assumption of the lower
bound for $m(v)$ or $m'(v)$.)

\begin{thm}\label{thm m}
Assume RH. There exists an absolute constant $c>0$, such that
$$m(v/2\pi)\le m'(v^2)$$ for all $v<c$\ .
\end{thm}

Theorem \ref{thm m} is a consequence of the following result.
\begin{thm}\label{thm box and m}
Assume RH. There exists an absolute constant $c>0$ such that for any
$v<c$ the following holds: For all large $\g$ with $\g^+-\g< v/\log
\g$, the box
$$
\{s=\sigma+it: \frac{1}{2}<\sigma<\frac{1}{2}+\frac{v^2}{4\log\g},
\g\le t\le \g^+\}
$$
contains exactly one zero of $\zeta'(s)$. Moreover, the zero is not
on the boundary of the box.
\end{thm}

Theorem \ref{thm box and m} can be viewed as a refinement of a
result of Zhang (see Theorem 3 in~\cite{Zha}). Roughly speaking,
Zhang's result states that (on RH) if there occurs a small gap of
consecutive zeros of $\zeta(s)$, then we may find a zero of
$\zeta'(s)$ nearby. The main difference between Zhang's result and
our Theorem \ref{thm box and m} is that, Theorem~\ref{thm box and m}
tells a more accurate location of $\rho'$: it allows us to put
$\rho'$ between consecutive zeros of $\zeta(s)$. Consequently,
different pairs of $(\g, \g^+)$ will give different $\rho'$ 's, and
this is how we shall deduce Theorem~\ref{thm m}.

Theorem \ref{thm box and m} should also be compared with a result of
Soundararajan (see Lemma \ref{lem sound at most 1 } in Section
\ref{sec preliminary}). Roughly speaking, Soundararajan's result
states that there is \emph{at most} one $\rho'$ in a certain
rectangular region between consecutive zeros of $\zeta(s)$. We will
prove a complementary result (see Lemma \ref{lem zero box} in
Section \ref{sec preliminary}), which asserts that there is \emph{at
least} one $\rho'$ in a certain rectangular region. Apart from
technical details, Theorem \ref{thm box and m} is essentially a
combination of these two results.

Our next result considers the unsolved part of Soundararajan's
Conjecture.
\begin{thm}\label{thm sound conj} Assume RH. For
$\beta'>1/2$ and $\g\le \g'<\g^+$, we have
$$
\g^+-\g\ll \sqrt{(\beta'-1/2)\log\g'}.
$$
In particular, we have
\begin{align}\label{eq weak sound conj}
\liminf_{\substack{\g'\rightarrow\infty\\\beta'>1/2}}(\beta'-1/2)(\log\g')^3=0\Longrightarrow
\liminf_{\g\rightarrow\infty}(\g^+-\g)\log\g=0.
\end{align}
\end{thm}

This should be compared with a result of M. Z. Garaev and C. Y.
Y{\i}ld{\i}r{\i}m~\cite{GYi}, which states that (on RH)
\begin{align}\label{eq G-Y weak sound conj}
\liminf_{\g'\rightarrow\infty}(\beta'-1/2)(\log\g')(\log\log\g')^2=0\Longrightarrow
\liminf_{\g\rightarrow\infty}(\g_{n+1}-\g_n)\log\g_n=0.
\end{align}
Note that the right-hand side of \eqref{eq G-Y weak sound conj}
would be trivially true if $\zeta(s)$ has infinitely many multiple
zeros. Our result \eqref{eq weak sound conj} overcomes this possible
issue since $\g^+>\g$ are ordinates of distinct zeros, but at the
cost that we require a stronger assumption in the left-hand side.
Note also that both \eqref{eq weak sound conj} and \eqref{eq G-Y
weak sound conj} are weaker than Soundararajan's Conjecture.

The paper is organized as follows. In Section \ref{sec preliminary}
we state preliminary results required in proving our main theorems.
We will then deduce Theorem \ref{thm box and
 m}, Theorem \ref{thm m} and Theorem \ref{thm sound conj} in order in Section \ref{sec proof thm},
 using results from Section \ref{sec preliminary}.
 The last two sections are devoted to the proofs of preliminary
 results.

\medskip

%%%%%%%%%%%%%%%%%%%%%%%%%%%%%%%%%%%%%%%%%%
%%%%%%%%%%%%%%%%%%%%%%%%%%%%%%%%%%%%%%%%%%
%%%%%%%%%%%%%%%%%%%%%%%%%%%%%%%%%%%%%%%%%%
%%%%%%%%%%%%%%%%%%%%%%%%%%%%%%%%%%%%%%%%%%
%%%%%%%%%%%%%%%%%%%%%%%%%%%%%%%%%%%%%%%%%%
%%%%%%%%%%%%%%%%%%%%%%%%%%%%%%%%%%%%%%%%%%
%%%%%%%%%%%%%%%%%%%%%%%%%%%%%%%%%%%%%%%%%%
%%%%%%%%%%%%%%%%%%%%%%%%%%%%%%%%%%%%%%%%%%
%%%%%%%%%%%%%%%%%%%%%%%%%%%%%%%%%%%%%%%%%%
%%%%%%%%%%%%%%%%%%%%%%%%%%%%%%%%%%%%%%%%%%
%%%%%%%%%%%%%%%%%%%%%%%%%%%%%%%%%%%%%%%%%%
%%%%%%%%%%%%%%%%%%%%%%%%%%%%%%%%%%%%%%%%%%
%%%%%%%%%%%%%%%%%%%%%%%%%%%%%%%%%%%%%%%%%%
%%%%%%%%%%%%%%%%%%%%%%%%%%%%%%%%%%%%%%%%%%
%%%%%%%%%%%%%%%%%%%%%%%%%%%%%%%%%%%%%%%%%%
%%%%%%%%%%%%%%%%%%%%%%%%%%%%%%%%%%%%%%%%%%
%%%%%%%%%%%%%%%%%%%%%%%%%%%%%%%%%%%%%%%%%%
%%%%%%%%%%%%%%%%%%%%%%%%%%%%%%%%%%%%%%%%%%
%%%%%%%%%%%%%%%%%%%%%%%%%%%%%%%%%%%%%%%%%%
%%%%%%%%%%%%%%%%%%%%%%%%%%%%%%%%%%%%%%%%%%
%%%%%%%%%%%%%%%%%%%%%%%%%%%%%%%%%%%%%%%%%%
%%%%%%%%%%%%%%%%%%%%%%%%%%%%%%%%%%%%%%%%%%
%%%%%%%%%%%%%%%%%%%%%%%%%%%%%%%%%%%%%%%%%%
%%%%%%%%%%%%%%%%%%%%%%%%%%%%%%%%%%%%%%%%%%
%%%%%%%%%%%%%%%%%%%%%%%%%%%%%%%%%%%%%%%%%%
%%%%%%%%%%%%%%%%%%%%%%%%%%%%%%%%%%%%%%%%%%
%%%%%%%%%%%%%%%%%%%%%%%%%%%%%%%%%%%%%%%%%%
%%%%%%%%%%%%%%%%%%%%%%%%%%%%%%%%%%%%%%%%%%
%%%%%%%%%%%%%%%%%%%%%%%%%%%%%%%%%%%%%%%%%%
%%%%%%%%%%%%%%%%%%%%%%%%%%%%%%%%%%%%%%%%%%
%%%%%%%%%%%%%%%%%%%%%%%%%%%%%%%%%%%%%%%%%%
%%%%%%%%%%%%%%%%%%%%%%%%%%%%%%%%%%%%%%%%%%
%%%%%%%%%%%%%%%%%%%%%%%%%%%%%%%%%%%%%%%%%%
%%%%%%%%%%%%%%%%%%%%%%%%%%%%%%%%%%%%%%%%%%
%%%%%%%%%%%%%%%%%%%%%%%%%%%%%%%%%%%%%%%%%%

\section{Preliminary results}\label{sec preliminary}

Let $\eta(s)=\pi^{-s/2}\Gamma(s/2)\zeta'(s)$, and define
\begin{equation}\label{def of
F} F(t) = \left\{
\begin{array}{ccc} & -\Re\ \dfrac{\eta'}{\eta}\big(1/2+it\big),  & \textrm{ if }\ \eta(1/2+it)\ne
0,\\ & \lim_{v\rightarrow t} F(v),  & \textrm{ otherwise.}
\end{array}\right.
\end{equation}

It is well-defined (see ~\cite{Zha}), and serves as a key role in
our proofs. The following five results were obtained by Zhang~\cite{Zha}. \\
(i) The limit in (\ref{def of F}) exists. Namely, $F(t)$ is
well-defined. Moreover, $F(t)$ is continuous.
\\
(ii) We have $F(t)=F_1(t)-F_2(t)+O(1)$, where
\begin{equation}\label{def F1}
F_1(t)= - \sum_{\beta'>1/2}\Re\frac{1}{1/2+it-\rho'}\ ,
\end{equation}
and
\begin{equation}\label{def F2}
F_2(t)= \sum_{0<\beta'<1/2}\Re\frac{1}{1/2+it-\rho'}\ .
\end{equation}
(iii) If $\rho=1/2+i\g$ is a simple zero of $\zeta(s)$ where $\g>0$,
then
 $F(\g)=\tfrac{1}{2}\log\g+O(1)$.\\
(iv) We have
\begin{align}\label{eq zhang int < pi}
\int_{\g}^{\g^+} F(t)dt \le \pi.
\end{align}
(v) If both $1/2+i\g$ and $1/2+i\g^+$ are simple zeros of
$\zeta(s)$, then
$$
\int_{\g}^{\g^+} F(t)dt \equiv 0\pmod{\pi}.
$$

We require some further results of $F(t)$ given by the following
four lemmas.
\begin{lem}\label{lem F} Let $\rho=1/2+i\g$ be a zero of $\zeta(s)$
with multiplicity $m(\rho)=m$. Then we have
$$
F(\g) = \frac{1}{2m} \log \g + O\big(\frac{1}{m}\big).
$$
\end{lem}

\begin{lem}\label{lem arg equiv}
Let $\rho=1/2+i\g$ be a zero of $\zeta(s)$. Then we have
$$
\lim_{\substack{v\rightarrow \g\\v>\g}}\arg\eta(1/2+iv)\equiv
\lim_{\substack{v\rightarrow \g\\v<\g}}\arg\eta(1/2+iv)\equiv \pi/2
\pmod{\pi}.
$$
\end{lem}

\begin{lem}\label{lem F int}
Assume RH. We have
$$
\int_{\g}^{\g^+} F(t)dt =\pi
$$
for all large $\g$.
\end{lem}

\begin{lem}\label{lem F F1}
 Assume RH, and let $F_1(t)$ be defined in
\eqref{def F1}. Then we have
$$
F(t)=F_1(t)+\log 2/2 +O(1/t).
$$
\end{lem}

For $\Re(\rho')> 1/2$, let $\theta(\rho', t_1, t_2)\in (0, \pi)$ be
the argument of the angle at $\rho'$ with rays through $1/2 + it_1$
and $1/2 + it_2$ respectively. The following lemma is also needed.

\begin{lem}\label{lem theta}
Assume RH. We have
$$
\sum_{1/2<\Re{\rho'}} \theta(\rho', \g, \g^+)
+\frac{\log2}{2}(\g^+-\g) + O\Big(\frac{\g^+-\g}{\g}\Big)=\pi
$$
for all large $\g$.
\end{lem}

Next we state two useful results of Soundararajan.

\begin{lem}\label{lem sound distance}
For $\beta'>1/2$ and $\g'>0$, we have
$$
\big|\rho-\rho'\big|\ge \sqrt{2(\beta'-1/2)/\log \g}.
$$
\end{lem}
See Lemma 2.1 in \cite{Sou}.

\begin{lem}\label{lem sound at most 1 } Assume RH. The box
$$
\{s=\sigma+it: 1/2<\sigma<1/2+1/\log \g, \g\le t\le \g^+\}
$$
contains at most one zero of $\zeta'(s)$.
\end{lem}
See Proposition 1.6 in \cite{Sou}.

The following result is also required, which can be viewed as a
partial complement to Lemma \ref{lem sound at most 1 }.

\begin{lem}\label{lem zero box}
Assume RH. Let $a$ be any constant less than $\pi/3$. There exists a
constant $\g_0(a)$ such that for $\g>\g_0(a)$ with $\g^+-\g<a/\log
\g$, the box
$$
\{s=\sigma+it: 1/2<\sigma<1/2+2.5(\g^+-\g), \g< t<\g^+\}
$$
contains a zero of $\zeta'(s)$.
\end{lem}

Combining the above two results we immediately obtain
\begin{cor}\label{cor box 0.4}
Assume RH. For large $\g$ with $\g^+-\g<0.4/\log \g$, the box
$$
\{s=\sigma+it: 1/2<\sigma<1/2+1/\log\g, \g\le t\le \g^+\}
$$
contains exactly one zero of $\zeta'(s)$. Moreover, the zero is not
on the boundary of the box.
\end{cor}

Lastly, it is interesting to record our final proposition, whose
proof we shall omit because of its similarity to that of Lemma
\ref{lem zero box}. (Also, we do not need it in proving our main
theorems.)
\begin{prop}\label{prop strip}
Assume RH. Let $\Delta$ be any constant less than $\pi$. There
exists a constant $\g_0(\Delta)$ such that for $\g>\g_0(\Delta)$
with $\g^+-\g<\Delta/\log \g$, the strip $\{s=\sigma+it:
\g<t<\g^+\}$ contains a zero of $\zeta'(s)$.
\end{prop}

%%%%%%%%%%%%%%%%%%%%%%%%%%%%%%%%%%%%%%%%%%
%%%%%%%%%%%%%%%%%%%%%%%%%%%%%%%%%%%%%%%%%%
%%%%%%%%%%%%%%%%%%%%%%%%%%%%%%%%%%%%%%%%%%
%%%%%%%%%%%%%%%%%%%%%%%%%%%%%%%%%%%%%%%%%%
%%%%%%%%%%%%%%%%%%%%%%%%%%%%%%%%%%%%%%%%%%
%%%%%%%%%%%%%%%%%%%%%%%%%%%%%%%%%%%%%%%%%%
%%%%%%%%%%%%%%%%%%%%%%%%%%%%%%%%%%%%%%%%%%
%%%%%%%%%%%%%%%%%%%%%%%%%%%%%%%%%%%%%%%%%%
%%%%%%%%%%%%%%%%%%%%%%%%%%%%%%%%%%%%%%%%%%
%%%%%%%%%%%%%%%%%%%%%%%%%%%%%%%%%%%%%%%%%%
%%%%%%%%%%%%%%%%%%%%%%%%%%%%%%%%%%%%%%%%%%
%%%%%%%%%%%%%%%%%%%%%%%%%%%%%%%%%%%%%%%%%%

%%%%%%%%%%%%%%%%%%%%%%%%%%%%%%%%%%%%%%%%%%
%%%%%%%%%%%%%%%%%%%%%%%%%%%%%%%%%%%%%%%%%%
%%%%%%%%%%%%%%%%%%%%%%%%%%%%%%%%%%%%%%%%%%
%%%%%%%%%%%%%%%%%%%%%%%%%%%%%%%%%%%%%%%%%%
%%%%%%%%%%%%%%%%%%%%%%%%%%%%%%%%%%%%%%%%%%
%%%%%%%%%%%%%%%%%%%%%%%%%%%%%%%%%%%%%%%%%%
%%%%%%%%%%%%%%%%%%%%%%%%%%%%%%%%%%%%%%%%%%
%%%%%%%%%%%%%%%%%%%%%%%%%%%%%%%%%%%%%%%%%%
%%%%%%%%%%%%%%%%%%%%%%%%%%%%%%%%%%%%%%%%%%
%%%%%%%%%%%%%%%%%%%%%%%%%%%%%%%%%%%%%%%%%%
%%%%%%%%%%%%%%%%%%%%%%%%%%%%%%%%%%%%%%%%%%
%%%%%%%%%%%%%%%%%%%%%%%%%%%%%%%%%%%%%%%%%%
%%%%%%%%%%%%%%%%%%%%%%%%%%%%%%%%%%%%%%%%%%
%%%%%%%%%%%%%%%%%%%%%%%%%%%%%%%%%%%%%%%%%%
%%%%%%%%%%%%%%%%%%%%%%%%%%%%%%%%%%%%%%%%%%
%%%%%%%%%%%%%%%%%%%%%%%%%%%%%%%%%%%%%%%%%%
%%%%%%%%%%%%%%%%%%%%%%%%%%%%%%%%%%%%%%%%%%
%%%%%%%%%%%%%%%%%%%%%%%%%%%%%%%%%%%%%%%%%%
%%%%%%%%%%%%%%%%%%%%%%%%%%%%%%%%%%%%%%%%%%
%%%%%%%%%%%%%%%%%%%%%%%%%%%%%%%%%%%%%%%%%%
%%%%%%%%%%%%%%%%%%%%%%%%%%%%%%%%%%%%%%%%%%
%%%%%%%%%%%%%%%%%%%%%%%%%%%%%%%%%%%%%%%%%%
%%%%%%%%%%%%%%%%%%%%%%%%%%%%%%%%%%%%%%%%%%
%%%%%%%%%%%%%%%%%%%%%%%%%%%%%%%%%%%%%%%%%%

\section{Proofs of Theorems \ref{thm box and
m}, \ref{thm m} and \ref{thm sound conj}}\label{sec proof thm}

\proof[Proof \,of\, Theorem \ref{thm box and m}] Let $v<0.4$ and
suppose $\g^+-\g<v/\log \g$. By Corollary \ref{cor box 0.4}, the box
$$
\{s=\sigma+it: 1/2<\sigma<1/2+1/\log\g, \g\le t\le \g^+\}
$$
contains exactly one zero of $\zeta'(s)$, say $\rho'=\beta'+i\g'$,
and it is not on the boundary. Thus, to prove the theorem, it
suffices to prove that
 \begin{align}\label{eq v^2}
 \beta'-1/2 < \frac{v^2}{4\log
\g}\ .\end{align}

Note that $\beta'-1/2<1/\log\g$. Without loss of generality, we may
assume $\g^+-\g'\ge\g'-\g$, and hence $\g^+-\g\ge 2(\g'-\g)$.

By Lemma \ref{lem sound distance}, we have
$$ \big|\rho'-\rho\big| \ge
\sqrt{\frac{2(\beta'-1/2)}{\log\g}},
$$
 and this gives
$$
(\g'-\g)^2\ge \frac{2(\beta'-1/2)}{\log\g} - (\beta'-1/2)^2\ge
\frac{\beta'-1/2}{\log\g},
$$
since $\beta'-1/2<1/\log\g$.

Hence, we get
$$
\frac{v^2}{\log^2\g}>(\g^+-\g)^2\ge 4(\g'-\g)^2\ge
\frac{4(\beta'-1/2)}{\log\g}.
$$

This gives \eqref{eq v^2} and completes the proof.

\qed

\medskip

%%%%%%%%%%%%%%%%%%%%%%%%%%%%%%%%%%%%%%%%%%
%%%%%%%%%%%%%%%%%%%%%%%%%%%%%%%%%%%%%%%%%%
%%%%%%%%%%%%%%%%%%%%%%%%%%%%%%%%%%%%%%%%%%
%%%%%%%%%%%%%%%%%%%%%%%%%%%%%%%%%%%%%%%%%%
%%%%%%%%%%%%%%%%%%%%%%%%%%%%%%%%%%%%%%%%%%
%%%%%%%%%%%%%%%%%%%%%%%%%%%%%%%%%%%%%%%%%%
%%%%%%%%%%%%%%%%%%%%%%%%%%%%%%%%%%%%%%%%%%
%%%%%%%%%%%%%%%%%%%%%%%%%%%%%%%%%%%%%%%%%%
%%%%%%%%%%%%%%%%%%%%%%%%%%%%%%%%%%%%%%%%%%
%%%%%%%%%%%%%%%%%%%%%%%%%%%%%%%%%%%%%%%%%%
%%%%%%%%%%%%%%%%%%%%%%%%%%%%%%%%%%%%%%%%%%
%%%%%%%%%%%%%%%%%%%%%%%%%%%%%%%%%%%%%%%%%%

\proof[Proof\, of\, Theorem \ref{thm m}]

Take $c$ to be the same as in Theorem \ref{thm box and m}, and let
$v<c$. Define
$$
\mathscr{S} = \{ n: T\le \g_n \le 2T, \g_{n+1}-\g_n\le \frac{v}{\log
T}\}
$$
and
$$
\mathscr{T} = \{ m: T\le \g'_m \le 2T, \beta'_m-\frac{1}{2}\le
\frac{v^2}{\log T}\}  .
$$

Recall that
$$
m'(v)=\liminf_{T\rightarrow \infty}\frac{2\pi}{T\log T}\sum_{T\le
\g'\le 2T} w\bigg(\frac{(\beta'-1/2)\log T}{v}\bigg)
$$
and
$$
m(v)=\liminf_{T\rightarrow \infty}\frac{2\pi}{T\log T}\sum_{T\le
\g_n\le 2T} w\bigg(\frac{(\g_{n+1}-\g_n)\log T}{2\pi v}\bigg),
$$
where $w(x)$ is the indicator function of $[0, 1]$.

Thus, to prove $m(v/2\pi)\le m'(v^2)$, it suffices to show that
$|\mathscr{S}|-1\le |\mathscr{T}|$ for all large $T$.

Write $\mathscr{S}=\mathscr{S}_1\cup\mathscr{S}_2$, where
$$
\mathscr{S}_1 =  \{ n: T\le \g_n \le 2T, \g_{n+1}-\g_n = 0\},
$$
and
$$
\mathscr{S}_2 = \{ n: T\le \g_n \le 2T, 0<\g_{n+1}-\g_n\le
\frac{v}{\log T}\}.
$$
Similarly, write $\mathscr{T}=\mathscr{T}_1\cup\mathscr{T}_2$, where
$$
\mathscr{T}_1 = \{ k: T\le \g'_k \le 2T, \beta'_k-\frac{1}{2}=0\},
$$
and
$$
\mathscr{T}_2 = \{ k: T\le \g'_k \le 2T, 0<\beta'_k-\frac{1}{2}\le
\frac{v^2}{\log T}\}.
$$
We clearly have $|\mathscr{S}|=|\mathscr{S}_1|+|\mathscr{S}_2|$ and
$|\mathscr{T}|=|\mathscr{T}_1|+|\mathscr{T}_2|$. Moreover, it is
easy to see that there is a bijection between $\mathscr{S}_1 $ and
$\mathscr{T}_1 $. Namely, we have $|\mathscr{S}_1 |=|\mathscr{T}_1
|$.

Let $\alpha$ be the largest element in $\mathscr{S}_2$. We show that
there is an injective mapping from $\mathscr{S}_2-\{\alpha\}$ to
$\mathscr{T}_2$. This will give $|\mathscr{S}_2 | -1
\le|\mathscr{T}_2 |$, and therefore completes the proof.

By the definition of $\mathscr{S}_2$, we have $\g_{n_1}> \g_{n_2}$
if $n_1, n_2\in \mathscr{S}_2$ and $n_1>n_2$. Thus, if $n\in
\mathscr{S}_2-\{\alpha\}$, then $\g_n < \g_{\alpha}$, which implies
that $\g_n^+\le \g_{\alpha}\le 2T$.

Take any $n\in \mathscr{S}_2-\{\alpha\}$. It follows from the
definition of $\mathscr{S}_2$ that $\g_n^+-\g_n\le v/\log T \le
2v/\log \g_n.$ Applying Theorem \ref{thm box and m}, we see that the
box
$$
\{s=\sigma+it: 1/2<\sigma<1/2+v^2/\log\g_n, \g_n< t <\g_n^+\}
$$
contains a zero of $\zeta'(s)$, say $\rho'_{k(n)}$. Since
$$
0<\beta'_{k(n)}-\frac{1}{2}<\frac{v^2}{\log\g_n}\le\frac{v^2}{\log
T}, \ \ \ \textrm{ and } \ \ T\le \g_n <\g'_{k(n)}<\g_n^+\le 2T,
$$
it follows that $k(n) \in \mathscr{T}_2$.

We take the mapping $\phi: \mathscr{S}_2-\{\alpha\} \longrightarrow
\mathscr{T}_2$ via $\phi(n)=k(n)$. It remains to show that $\phi$ is
injective. Suppose this is not the case, then we would have
$k(n_1)=k(n_2)$ for some $n_1 > n_2$. But that would imply
$$\g'_{k(n_1)}>\g_{n_1}\ge \g_{n_2}^+>\g'_{k(n_2)}=\g'_{k(n_1)},$$ a
contradiction. This completes the proof.

\qed

\medskip

%%%%%%%%%%%%%%%%%%%%%%%%%%%%%%%%%%%%%%%%%%
%%%%%%%%%%%%%%%%%%%%%%%%%%%%%%%%%%%%%%%%%%
%%%%%%%%%%%%%%%%%%%%%%%%%%%%%%%%%%%%%%%%%%
%%%%%%%%%%%%%%%%%%%%%%%%%%%%%%%%%%%%%%%%%%
%%%%%%%%%%%%%%%%%%%%%%%%%%%%%%%%%%%%%%%%%%
%%%%%%%%%%%%%%%%%%%%%%%%%%%%%%%%%%%%%%%%%%
%%%%%%%%%%%%%%%%%%%%%%%%%%%%%%%%%%%%%%%%%%
%%%%%%%%%%%%%%%%%%%%%%%%%%%%%%%%%%%%%%%%%%
%%%%%%%%%%%%%%%%%%%%%%%%%%%%%%%%%%%%%%%%%%
%%%%%%%%%%%%%%%%%%%%%%%%%%%%%%%%%%%%%%%%%%
%%%%%%%%%%%%%%%%%%%%%%%%%%%%%%%%%%%%%%%%%%
%%%%%%%%%%%%%%%%%%%%%%%%%%%%%%%%%%%%%%%%%%

\proof[Proof \,of \,Theorem \ref{thm sound conj}] Let
$\rho'=\beta'+i\g'$ be a zero of $\zeta'(s)$ with $\beta'>1/2$ and
$\g'$ large. Let $\g$ and $\g^+$ be such that $\g<\g'<\g^+$.

By Lemma \ref{lem theta} we have
$$
\sum_{1/2<\Re{\lambda'}} \theta(\lambda', \g, \g^+)
+\frac{\log2}{2}(\g^+-\g) + O(\frac{\g^+-\g}{\g})=\pi,
$$
where the sum is over all zeros $\lambda'$ of $\zeta'(s)$ with real
part greater than 1/2. In particular, this gives
$$
\theta(\rho', \g, \g^+) <\pi -\frac{\log2}{2}(\g^+-\g) +
O(\frac{\g^+-\g}{\g}).
$$
Thus, it follows that
$$
\theta_1+\theta_2>\frac{\log2}{2}(\g^+-\g) +
O(\frac{\g^+-\g}{\g})\gg \g^+-\g,
$$
where $\theta_1$ and $\theta_2$ are the angles of at $\rho$ and
$\rho^+$, respectively, of the triangle $(\rho, \rho', \rho^+)$.

By Lemma \ref{lem sound distance}
$$
\big|\rho-\rho'\big|\gg\sqrt{(\beta'-1/2)/\log \g}.
$$
Therefore, we see that
$$
\sin\theta_1 =
\frac{\beta'-1/2}{\big|\rho-\rho'\big|}\ll\sqrt{(\beta'-1/2)\log
\g},
$$
and this gives
$$
\theta_1\ll\sqrt{(\beta'-1/2)\log \g}.
$$
Similarly, we have
$$ \theta_2\ll\sqrt{(\beta'-1/2)\log \g}.
$$
Thus, we see that
$$
\g^+-\g\ll \theta_1+\theta_2\ll\sqrt{(\beta'-1/2)\log \g}.
$$
This proves the theorem. \qed

\medskip
%%%%%%%%%%%%%%%%%%%%%%%%%%%%%%%%%%%%%%%%%%
%%%%%%%%%%%%%%%%%%%%%%%%%%%%%%%%%%%%%%%%%%
%%%%%%%%%%%%%%%%%%%%%%%%%%%%%%%%%%%%%%%%%%
%%%%%%%%%%%%%%%%%%%%%%%%%%%%%%%%%%%%%%%%%%
%%%%%%%%%%%%%%%%%%%%%%%%%%%%%%%%%%%%%%%%%%
%%%%%%%%%%%%%%%%%%%%%%%%%%%%%%%%%%%%%%%%%%
%%%%%%%%%%%%%%%%%%%%%%%%%%%%%%%%%%%%%%%%%%
%%%%%%%%%%%%%%%%%%%%%%%%%%%%%%%%%%%%%%%%%%
%%%%%%%%%%%%%%%%%%%%%%%%%%%%%%%%%%%%%%%%%%
%%%%%%%%%%%%%%%%%%%%%%%%%%%%%%%%%%%%%%%%%%
%%%%%%%%%%%%%%%%%%%%%%%%%%%%%%%%%%%%%%%%%%
%%%%%%%%%%%%%%%%%%%%%%%%%%%%%%%%%%%%%%%%%%
%%%%%%%%%%%%%%%%%%%%%%%%%%%%%%%%%%%%%%%%%%
%%%%%%%%%%%%%%%%%%%%%%%%%%%%%%%%%%%%%%%%%%
%%%%%%%%%%%%%%%%%%%%%%%%%%%%%%%%%%%%%%%%%%
%%%%%%%%%%%%%%%%%%%%%%%%%%%%%%%%%%%%%%%%%%
%%%%%%%%%%%%%%%%%%%%%%%%%%%%%%%%%%%%%%%%%%
%%%%%%%%%%%%%%%%%%%%%%%%%%%%%%%%%%%%%%%%%%
%%%%%%%%%%%%%%%%%%%%%%%%%%%%%%%%%%%%%%%%%%
%%%%%%%%%%%%%%%%%%%%%%%%%%%%%%%%%%%%%%%%%%
%%%%%%%%%%%%%%%%%%%%%%%%%%%%%%%%%%%%%%%%%%
%%%%%%%%%%%%%%%%%%%%%%%%%%%%%%%%%%%%%%%%%%
%%%%%%%%%%%%%%%%%%%%%%%%%%%%%%%%%%%%%%%%%%
%%%%%%%%%%%%%%%%%%%%%%%%%%%%%%%%%%%%%%%%%%

\section{Proofs of Lemmas \ref{lem F}, \ref{lem arg equiv}, \ref{lem F int}, \ref{lem F F1} and \ref{lem theta}} \label{sec proof
lem}

\proof[Proof \,of\, Lemma \ref{lem F}]

Following Zhang \cite{Zha}, we let $\xi(s)=h(s)\zeta(s)$ and
$\eta(s)=h(s)\zeta'(s)$, where $h(s)=\pi^{-s/2}\Gamma(s/2)$. Note
that for any integer $n>0$, we have
$$
i^n\xi^{(n)}(1/2+it)\in \mathbb{R}
$$
by the functional equation.

Suppose that $\rho=1/2+i\g$ is a zero of $\zeta(s)$ with
multiplicity $m$. Then we have
$$
\zeta(\rho)=\zeta'(\rho)=\cdots=\zeta^{(m-1)}(\rho)=0, \ \zeta^{(m)}(\rho)\ne 0.
$$
It follows from Leibniz' law that
$$
\xi(\rho)=\xi'(\rho)=\cdots=\xi^{(m-1)}(\rho)=\eta(\rho)=\eta'(\rho)=\cdots=\eta^{(m-2)}(\rho)=0, \ \xi^{(m)}(\rho)=\eta^{(m-1)}(\rho)\ne 0.
$$
In particular, we see that $i^m\eta^{(m-1)}(\rho)\in \mathbb{R}$.

Write $\eta(x+iy)=\eta(x,y)$. If $\eta(s)$ is holomorphic at
$s=x+iy$, then we have
\begin{align} \label{eta eta}
\eta'(x+iy) & = \eta_y(x,y)\cdot(i)^{-1},\nonumber\\
\eta^{(2)}(x+iy) & = \eta_{yy}(x,y)\cdot(i)^{-2},\nonumber\\
& \dots \nonumber\\
\eta^{(k)}(x+iy) & = \eta_{y^k}(x,y)\cdot(i)^{-k}, \ \ \forall k\in
\mathbb{Z}^+.
\end{align}
Now since $\eta(\rho)=\eta'(\rho)=\cdots=\eta^{(m-2)}(\rho)=0$ and
$\eta^{(m-1)}(\rho)\ne 0$, it follows that
$$
\eta(1/2, \g)=\eta_y(1/2, \g)=\cdots=\eta_{y^{m-2}}(1/2, \g)=0, \ \
\eta_{y^{m-1}}(1/2, \g)\ne 0.
$$
Therefore, when $t$ is in a neighborhood of $\g$ in which $\eta(s)$ is holomorphic, we can write
\begin{equation}\label{eta p}
\eta(1/2, t)=(t-\g)^{m-1}\cdot p(t)
\end{equation}
 for some function $p$ with
$p(\g)\ne 0$. Hence, we have
\begin{align*}
\eta'(1/2+it) & =   \eta_y(1/2,t)\cdot i^{-1}    \\
&= \bigg((t-\g)^{m-1}p'(t)+(m-1)(t-\g)^{m-2}p(t)\bigg) \cdot i^{-1}.
\end{align*}

This gives
$$
\frac{\eta'}{\eta}(1/2+it)=i^{-1}\frac{p'}{p}(t)+i^{-1}(m-1)(t-\g)^{-1},
$$
and in particular,
$$
-\Re\ \frac{\eta'}{\eta}(1/2+it)=-\Im\ \frac{p'}{p}(t).
$$

From this and the definition of $F(t)$, we see that
$$
 F(\g)= -\Im\ \frac{p'}{p}(\g)
 $$
 since $p(\g)\ne 0$.

Now Leibniz' law gives us
$$
\xi^{(m+1)}(\rho)=h\zeta^{(m+1)}(\rho)+(m+1)h'\zeta^{(m)}(\rho)
$$
and
$$
\xi^{(m)}(\rho)=h\zeta^{(m)}(\rho).
$$
It follows that
$$
\frac{\xi^{(m+1)}}{\xi^{(m)}(\rho)}=\frac{\zeta^{(m+1)}}{\zeta^{(m)}(\rho)}+(m+1)\frac{h'}{h}(\rho).
$$
Similarly, we obtain
$$
\frac{\eta^{(m)}}{\eta^{(m-1)}(\rho)}=\frac{\zeta^{(m+1)}}{\zeta^{(m)}(\rho)}+m\frac{h'}{h}(\rho).
$$
This gives
$$
\frac{\xi^{(m+1)}}{\xi^{(m)}(\rho)}- \frac{\eta^{(m)}}{\eta^{(m-1)}(\rho)}= \frac{h'}{h}(\rho).
$$

On the other hand, by (\ref{eta eta}) and  (\ref{eta p}) we can
easily compute that
$$
\frac{\eta^{(m)}}{\eta^{(m-1)}(\rho)}=i^{-1}m\frac{p'}{p}(\g).
$$
Combining the above two formulas, and taking the real part, we see
that
\begin{align*}
F(\g)=\frac{1}{m}\Re\frac{h'}{h}(\rho) =\frac{1}{2m} \log \g +
O(1/m),
\end{align*}
where we apply Stirling's formula in the last equality. This
completes the proof of Lemma \ref{lem F}.  \qed

\medskip

%%%%%%%%%%%%%%%%%%%%%%%%%%%%%%%%%%%%%%%%%%
%%%%%%%%%%%%%%%%%%%%%%%%%%%%%%%%%%%%%%%%%%
%%%%%%%%%%%%%%%%%%%%%%%%%%%%%%%%%%%%%%%%%%
%%%%%%%%%%%%%%%%%%%%%%%%%%%%%%%%%%%%%%%%%%
%%%%%%%%%%%%%%%%%%%%%%%%%%%%%%%%%%%%%%%%%%
%%%%%%%%%%%%%%%%%%%%%%%%%%%%%%%%%%%%%%%%%%
%%%%%%%%%%%%%%%%%%%%%%%%%%%%%%%%%%%%%%%%%%
%%%%%%%%%%%%%%%%%%%%%%%%%%%%%%%%%%%%%%%%%%
%%%%%%%%%%%%%%%%%%%%%%%%%%%%%%%%%%%%%%%%%%
%%%%%%%%%%%%%%%%%%%%%%%%%%%%%%%%%%%%%%%%%%
%%%%%%%%%%%%%%%%%%%%%%%%%%%%%%%%%%%%%%%%%%
%%%%%%%%%%%%%%%%%%%%%%%%%%%%%%%%%%%%%%%%%%
%%%%%%%%%%%%%%%%%%%%%%%%%%%%%%%%%%%%%%%%%%
%%%%%%%%%%%%%%%%%%%%%%%%%%%%%%%%%%%%%%%%%%
%%%%%%%%%%%%%%%%%%%%%%%%%%%%%%%%%%%%%%%%%%
%%%%%%%%%%%%%%%%%%%%%%%%%%%%%%%%%%%%%%%%%%

\proof[Proof\, of \,Lemma \ref{lem arg equiv}]

The case that $\rho=1/2+i\g$ is a simple zero of $\zeta(s)$ has been
treated in \cite{Zha}. Below we assume that  $\rho=1/2+i\g$  is a
multiple zero of $\zeta(s)$. Namely, we assume $\eta(1/2+i\g)=0$.

We use a temporary notation $\lim$ to denote
$\lim_{\substack{v\rightarrow \g\\v>\g}}$ (or
$\lim_{\substack{v\rightarrow \g\\v<\g}}$). It is easy to see that
\begin{align*}
\lim \arg \eta(1/2+it_1) & \equiv \lim\arctan\frac{\Im\ \eta(1/2+it_1)}{\Re\ \eta(1/2+it_1)}
\pmod{\pi}\\
& \equiv \arctan\lim \frac{\Im\ \eta(1/2+it_1)}{\Re\ \eta(1/2+it_1)} \pmod{\pi}.
\end{align*}

By \eqref{eta p} we have
$$
\Re\ \eta(1/2, t_1)=(t_1-\g)^{m-1}\cdot
\Re\ p(t_1)
$$
and
$$
\Im\ \eta(1/2, t_1)=(t_1-\g)^{m-1}\cdot \Im\ p(t_1).
$$
It follows that
\begin{align*}
\lim \arg \eta(1/2+it_1) & \equiv \arctan\lim \frac{\Im\ \eta(1/2+it_1)}{\Re\ \eta(1/2+it_1)} \pmod{\pi}\\
& \equiv \arctan\lim\frac{(t_1-\g)^{m-1}\cdot \Im\ p(t_1)}{(t_1-\g)^{m-1}\cdot \Re\ p(t_1)} \pmod{\pi}\\
& \equiv  \arctan\lim\frac{ \Im\ p(t_1)}{\Re\ p(t_1)}  \pmod{\pi}.
\end{align*}

Recall that $i^{m}\eta^{(m-1)}(\rho)\in \mathbb{R}$ and that
$\eta_{y^{m-1}}(1/2, \g)\cdot (i)^{-(m-1)}= \eta^{(m-1)}(\rho)$.
Thus, we have
$$
i\cdot\eta_{y^{m-1}}(1/2, \g)\in \mathbb{R}.
$$
Since $\eta_{y^{m-1}}(1/2, \g)=(m-1)!\ p(\g)$, we have
$ip(\g)\in\mathbb{R}$, namely, $\Re\ p(\g)=0$. From this and the
fact that $p(\g)\ne 0$, we see that $\Im \ p(\g)\ne 0$. Therefore,
we have
\begin{align*}
\lim \arg \eta(1/2+it_1) & \equiv  \arctan\lim\frac{ \Im\ p(t_1)}{\Re\ p(t_1)} \pmod{\pi}\\
&  \equiv \pi/2 \pmod{\pi}.
\end{align*}
\qed

\medskip

%%%%%%%%%%%%%%%%%%%%%%%%%%%%%%%%%%%%%%%%%%
%%%%%%%%%%%%%%%%%%%%%%%%%%%%%%%%%%%%%%%%%%
%%%%%%%%%%%%%%%%%%%%%%%%%%%%%%%%%%%%%%%%%%
%%%%%%%%%%%%%%%%%%%%%%%%%%%%%%%%%%%%%%%%%%
%%%%%%%%%%%%%%%%%%%%%%%%%%%%%%%%%%%%%%%%%%
%%%%%%%%%%%%%%%%%%%%%%%%%%%%%%%%%%%%%%%%%%
%%%%%%%%%%%%%%%%%%%%%%%%%%%%%%%%%%%%%%%%%%
%%%%%%%%%%%%%%%%%%%%%%%%%%%%%%%%%%%%%%%%%%
%%%%%%%%%%%%%%%%%%%%%%%%%%%%%%%%%%%%%%%%%%
%%%%%%%%%%%%%%%%%%%%%%%%%%%%%%%%%%%%%%%%%%
%%%%%%%%%%%%%%%%%%%%%%%%%%%%%%%%%%%%%%%%%%
%%%%%%%%%%%%%%%%%%%%%%%%%%%%%%%%%%%%%%%%%%
%%%%%%%%%%%%%%%%%%%%%%%%%%%%%%%%%%%%%%%%%%
%%%%%%%%%%%%%%%%%%%%%%%%%%%%%%%%%%%%%%%%%%
%%%%%%%%%%%%%%%%%%%%%%%%%%%%%%%%%%%%%%%%%%
%%%%%%%%%%%%%%%%%%%%%%%%%%%%%%%%%%%%%%%%%%

\proof[Proof\, of\, Lemma \ref{lem F int}]

 By Zhang's result (i), we
have
\begin{equation}\label{eq int F}
\int_\g^{\g^+} F(t)dt=\lim\int_{t_1}^{t_2} F(t)dt ,
\end{equation}
where the limit is taken for $t_1>\g, t_1\rightarrow \g$ and
$t_2<\g^+, t_2\rightarrow \g^+$. By equation (2.16) in \cite{Zha}
\begin{equation}\label{eq int F appr}
\int_{t_1}^{t_2} F(t)dt = \arg \eta(1/2+it_1)- \arg \eta(1/2+it_2).
\end{equation}
It follows that
$$
\int_\g^{\g^+} F(t)dt=\lim \arg \eta(1/2+it_1)- \lim \arg
\eta(1/2+it_2).
$$
By Lemma \ref{lem arg equiv}, this is
\begin{align}\label{eq int equiv}
\int_\g^{\g^+} F(t)dt \equiv 0 \pmod{\pi}.
\end{align}

By Zhang's result (ii),
$$
F(t) = F_1(t) -F_2(t) + O(1).
$$
This comes from considering the Hadamard factorization for
$\eta(s)$. With a little more care we may actually get
$$
F(t) = F_1(t) -F_2(t) - C +O(1/t),
$$
for some constant $C$, and on RH this is
\begin{align}\label{eq F}
F(t) = F_1(t) - C +O(1/t).
\end{align}
The expression for $C$ is given by
\begin{align}\label{def C}
C = \Re \sum_{\beta'>0} \frac{1}{\rho'} +
\sum_{n=1}^{\infty}\bigg(\frac{1}{\rho'_n}+\frac{1}{2n}\bigg)-\frac{\log
\pi}{2}-\frac{C_0}{2}-2+\frac{\zeta''}{\zeta'}(0),
\end{align}
where $C_0$ is Euler's constant, and $\rho'_n\in(-2n-2, -2n)$ is a
real zero of $\zeta'(s)$ (see Theorem 9 in \cite{LMo}).

One can easily show that $C<0$ (unconditionally). In fact, by
equation (4) in \cite{Yil}, we have
$$
\Re\sum_{\beta'>0}\frac{1}{\rho'}<0.185.
$$
Also, it is easy to calculate that
$$
\sum_{n=1}^{\infty}\bigg(\frac{1}{2n}+ \frac{1}{\rho'_n}\bigg)<0.32.
$$
Inserting the values for other constants in \eqref{def C}, we get
$C<-0.17$.

By \eqref{eq F} we see that $F(t)>F_1(t)\ge 0$ for large $t$. It
follows that
$$
\int_\g^{\g^+}F(t)dt >0.
$$

Combining this with \eqref{eq zhang int < pi} and \eqref{eq int
equiv}, we obtain
$$
\int_\g^{\g^+}F(t)dt =\pi
$$
for all large $\g$. \qed

\medskip

%%%%%%%%%%%%%%%%%%%%%%%%%%%%%%%%%%%%%%%%%%
%%%%%%%%%%%%%%%%%%%%%%%%%%%%%%%%%%%%%%%%%%
%%%%%%%%%%%%%%%%%%%%%%%%%%%%%%%%%%%%%%%%%%
%%%%%%%%%%%%%%%%%%%%%%%%%%%%%%%%%%%%%%%%%%
%%%%%%%%%%%%%%%%%%%%%%%%%%%%%%%%%%%%%%%%%%
%%%%%%%%%%%%%%%%%%%%%%%%%%%%%%%%%%%%%%%%%%
%%%%%%%%%%%%%%%%%%%%%%%%%%%%%%%%%%%%%%%%%%
%%%%%%%%%%%%%%%%%%%%%%%%%%%%%%%%%%%%%%%%%%
%%%%%%%%%%%%%%%%%%%%%%%%%%%%%%%%%%%%%%%%%%
%%%%%%%%%%%%%%%%%%%%%%%%%%%%%%%%%%%%%%%%%%
%%%%%%%%%%%%%%%%%%%%%%%%%%%%%%%%%%%%%%%%%%
%%%%%%%%%%%%%%%%%%%%%%%%%%%%%%%%%%%%%%%%%%
%%%%%%%%%%%%%%%%%%%%%%%%%%%%%%%%%%%%%%%%%%
%%%%%%%%%%%%%%%%%%%%%%%%%%%%%%%%%%%%%%%%%%
%%%%%%%%%%%%%%%%%%%%%%%%%%%%%%%%%%%%%%%%%%
%%%%%%%%%%%%%%%%%%%%%%%%%%%%%%%%%%%%%%%%%%

\proof[Proof\, of\, Lemma \ref{lem F F1}]

Consider the integral
$$
\int_T^{2T} F(t)dt.
$$

Ordering the ordinates $T\le g_1<g_2<\dots<g_{N_d}\le 2T$ of
distinct zeros of $\zeta(s)$ in $[T, 2T]$, we can write the above
integral as
\begin{align}\label{eq int F split}
\int_T^{g_1}+\int_{g_1}^{g_2}+\cdots+\int_{g_{N_{d}-1}}^{g_{N_{d}}}+\int_{g_{N_{d}}}^{2T}
F(t)dt=\sum_{j=1}^{N_{d}}\int_{g_j}^{g_{j+1}}F(t)dt +O(1).
\end{align}
On the other hand, by \eqref{eq F} we have
$$
\int_T^{2T} F(t)dt=\int_T^{2T} F_1(t)dt-CT+O(1).
$$
By the definition of $F_1(t)$,
$$
\int_T^{2T}
F_1(t)dt=-\Re\sum_{\beta'>1/2}\int_T^{2T}\frac{1}{1/2+it-\rho'}dt=
\sum_{\beta'>1/2}\theta(\rho', T, 2T).
$$

Hence, we obtain
$$
\int_T^{2T} F(t)dt =   \sum_{1/2<\beta'}\theta(\rho', T, 2T)   -CT
+O(1).
$$

By standard estimates, one easily shows that
$$
 \sum_{1/2<\beta'}\theta(\rho', T, 2T)=
 \pi N_{d} - T\log2/2+O(T/\log T).
 $$
It follows that
$$
\int_T^{2T} F(t)dt =   \pi N_{d} - T\log2/2 -CT + O(T/\log T).
$$

This together with (\ref{eq int F split}) gives
$$
\sum_{j=1}^{N_{d}}\bigg( \pi - \int_{g_j}^{g_{j+1}}
F(t)dt\bigg)=\bigg(C+\frac{\log 2}{2}\bigg)T+O(T/\log T).
$$
By Lemma \ref{lem F int}, the left-hand side is 0 for all large $T$.
Thus we have $C=-\log 2/2$. The result now follows from \eqref{eq
F}.
 \qed

\medskip

%%%%%%%%%%%%%%%%%%%%%%%%%%%%%%%%%%%%%%%%%%
%%%%%%%%%%%%%%%%%%%%%%%%%%%%%%%%%%%%%%%%%%
%%%%%%%%%%%%%%%%%%%%%%%%%%%%%%%%%%%%%%%%%%
%%%%%%%%%%%%%%%%%%%%%%%%%%%%%%%%%%%%%%%%%%
%%%%%%%%%%%%%%%%%%%%%%%%%%%%%%%%%%%%%%%%%%
%%%%%%%%%%%%%%%%%%%%%%%%%%%%%%%%%%%%%%%%%%
%%%%%%%%%%%%%%%%%%%%%%%%%%%%%%%%%%%%%%%%%%
%%%%%%%%%%%%%%%%%%%%%%%%%%%%%%%%%%%%%%%%%%
%%%%%%%%%%%%%%%%%%%%%%%%%%%%%%%%%%%%%%%%%%
%%%%%%%%%%%%%%%%%%%%%%%%%%%%%%%%%%%%%%%%%%
%%%%%%%%%%%%%%%%%%%%%%%%%%%%%%%%%%%%%%%%%%
%%%%%%%%%%%%%%%%%%%%%%%%%%%%%%%%%%%%%%%%%%
%%%%%%%%%%%%%%%%%%%%%%%%%%%%%%%%%%%%%%%%%%
%%%%%%%%%%%%%%%%%%%%%%%%%%%%%%%%%%%%%%%%%%
%%%%%%%%%%%%%%%%%%%%%%%%%%%%%%%%%%%%%%%%%%
%%%%%%%%%%%%%%%%%%%%%%%%%%%%%%%%%%%%%%%%%%

\proof[Proof\, of\, Lemma \ref{lem theta}]

The result follows immediately from Lemma \ref{lem F int} and Lemma
\ref{lem F F1}. \qed

\section{Proof of Lemma \ref{lem zero box}}

\proof[Proof\, of\, Lemma \ref{lem zero box}]

Let $\g$ be large. Let $\mathscr O=1/2+i(\g+\g^+)/2$. That is,
$\mathscr O$ is the middle point of $\rho$ and $\rho^+$. Write
$\g^+-\g=2l$ and $d=5l$. Also write
\begin{align}\label{eq F1 sep}
F_1(t)=F_{11}(t)+F_{12}(t),
\end{align}
where
$$
F_{11}(t)=-\sum_{\substack{\beta'>1/2\\|\rho'-\mathscr
O|<d}}\Re\frac{1}{1/2+it-\rho'}=
\sum_{\substack{\beta'>1/2\\|\rho'-\mathscr
O|<d}}\frac{\beta'-1/2}{(\beta'-1/2)^2+(\g'-t)^2}\ ,
$$
and
$$
F_{12}(t)=- \sum_{\substack{\beta'>1/2\\|\rho'-\mathscr O|\ge
d}}\Re\frac{1}{1/2+it-\rho'}=\sum_{\substack{\beta'>1/2\\|\rho'-\mathscr
O|\ge d}}\frac{\beta'-1/2}{(\beta'-1/2)^2+(\g'-t)^2}\ .
$$
For $|\rho'-\mathscr O|\ge d$ and $\g\le t\le \g^+$, we have
$$
|\rho' - (1/2+it)| \ge |\rho' - \rho| - |\rho - (1/2+it)| \ge |\rho'
- \rho|/2
$$
and
$$
|\rho' - (1/2+it)| \le |\rho' - \rho| + |\rho - (1/2+it)| \le
3|\rho' - \rho|/2.
$$
It follows that
$$
4F_{12}(\g)/9\le F_{12}(t) \le 4F_{12}(\g).
$$
In particular, this gives
\begin{align}\label{eq bound F12}
F_{12}(t)\le 2\log\g +O(1)
\end{align}
in view of Lemma \ref{lem F}.

Now suppose that there is no zero of $\zeta'(s)$ in the box
\begin{align}\label{assmptn}
\{s=\sigma+it: 1/2<\sigma<1/2+2.5(\g^+-\g), \g< t< \g^+\}.
\end{align}
Then we may write
$$
F_{11}(t)=f(t)+g(t),
$$
where
$$
f(t)= \sum_{\substack{\beta'>1/2\\|\rho'-\mathscr O|<d\\\g'\ge
\g^+}}\frac{\beta'-1/2}{(\beta'-1/2)^2+(\g'-t)^2}\ ,$$ and
$$
g(t)= \sum_{\substack{\beta'>1/2\\|\rho'-\mathscr
 O|<d\\\g'\le\g}}\frac{\beta'-1/2}{(\beta'-1/2)^2+(\g'-t)^2}\ .$$
Observe that for $\g\le t\le\g^+$ we have $f(t)\le f(\g^+)$ and
$g(t)\le g(\g)$. It follows that
$$
F_{11}(t)\le  f(\g^+)+ g(\g)\le F_1(\g^+) + F_1(\g). \le \log \g
+O(1).
$$
Applying Lemma \ref{lem F}, we obtain
\begin{align}\label{eq bound F11}
F_{11}(t)\le \log \g +O(1).
\end{align}

Combining \eqref{eq F1 sep}, \eqref{eq bound F11} and \eqref{eq
bound F12}, we get
$$
F_1(t)\le 3\log \g +O(1)
$$
for $t\in[\g,\g^+]$.
It follows that
$$
\int_{\g}^{\g^+} F_1(t)dt \le 3(\g^+-\g)\log\g +O(\g^+-\g).
$$

On the other hand, by Lemma \ref{lem F int} we have
$$
\int_{\g}^{\g^+} F(t)dt =\pi,
$$
and by Lemma \ref{lem F F1} this is
$$
\int_{\g}^{\g^+} F_1(t)dt =\pi + O(\g^+-\g).
$$

Hence, we obtain
$$
\pi\le 3(\g^+-\g)\log\g +O(\g^+-\g).
$$
But since we are assuming $\g^+-\g<a/\log\g$ (with $a<\pi/3$ \ a
constant), the above inequality would give
$$
\pi\le (3 + O(\log^{-1} \g)) a.
$$
This can not hold if $\g$ is large enough (depending on $a$). Hence
the assumption (\ref{assmptn}) is false. This completes our proof.
\qed

\end{document}